\documentclass[11pt]{amsproc}
\usepackage{amsfonts,amssymb}
\usepackage{eucal}
\usepackage[all]{xy}
\begin{document}
\xyoption{all}
\newcommand{\C}{{\mathbb C}}
\newcommand{\N}{{\mathcal N}}
\newcommand{\R}{{\mathbb R}}
\newcommand{\T}{{\mathbb T}}
\newcommand{\Z}{{\mathbb Z}}
\newcommand{\res}{{\rm res}}

\title {Cobordism of involutions revisited, revisited}\bigskip
\author{Jack Morava}
\address{Department of Mathematics, Johns Hopkins University,
Baltimore, Maryland 21218}
\email{jack@math.jhu.edu}
\thanks{The author was supported in part by the NSF}
\subjclass{Primary 55-03, Secondary 55N22, 55N91}
\date {7 March 1998}
\begin{abstract} Boardman's work [3,4] on Conner and Floyd's five-halves
conjecture looks remarkably contemporary in retrospect. This note
reexamines some of that work from a perspective proposed recently
by Greenlees and Kriz. \end{abstract}

\maketitle

\section {Unoriented involutions}

\noindent
A large part of Boardman's argument can be summarized as a
commutative diagram 
$$\xymatrix{
0 \ar[r]& \N_{\Z/2\Z*}^{geo} \ar[r] \ar[d]& \oplus_{i \geq 0} 
\N_{i}(BO_{*-i}) \ar[r]^{\partial} \ar[d]^{J} & \widetilde \N_{*-1}(B\Z/2\Z) 
\ar[r] \ar[d] & 0 \\
0 \ar[r]& \N_{*}[[w]] \ar[r]& \N_{*}((w)) \ar[r]& \N_{*-1}[w^{-1}]
\ar[r]& 0 \; .}
$$ 
of short exact sequences: here $\N_{*}((w))$ is the ring of homogeneous 
formal Laurent  series over the unoriented cobordism ring, in an indeterminate 
$w$ of cohomological degree one; such series are permitted only finitely many 
negative powers of $w$. The lower left-hand arrow is inclusion of the 
ring of formal power series, and the lower right-hand arrow is the quotient
homomorphism. The maps across the top are geometric: the left-hand arrow
sends the class of an unoriented manifold with involution to its
class in the relative bordism of unoriented $\Z/2\Z$-manifolds-with-boundary 
with action free on the boundary, while the right-hand arrow $\partial$ sends 
such a manifold to its boundary, regarded as a manifold with involution 
classified by a map to $B\Z/2\Z$. It is a theorem of Conner and Floyd 
[7 \S 22] that such a relative manifold with involution is cobordant to 
the normal ball bundle of its fixed-point set, so the group in the top row 
center is the sum of bordism groups of classifying spaces for orthogonal 
groups, as indicated. 

\newpage

\noindent{\sc Remarks:}\medskip

i) Thom showed that unoriented cobordism is a generalized
Eilenberg-Mac Lane spectrum, but this is not true of equivariant
unoriented bordism, cf.\ [8]; \medskip

ii) the homomorphism $\partial$ is {\bf not} a derivation;
\medskip

iii) the right-hand vertical arrow is an isomorphism, while
the left-hand vertical arrow sends a closed manifold with
involution $V$ to the bordism class of its Borel construction,
regarded as an element of $\N^{-*}(B\Z/2\Z)$, cf.\ [13]; \medskip

iv) the middle vertical homomorphism $J$ {\bf is} a ring homomorphism; its
construction [4, Theorem 9] is one of Boardman's innovations. He concludes
that it is in fact a monomorphism [4, Corollary 17], so
$\N_{\Z/2\Z*}^{geo}$ might be described as the subring of
classes $[V]$ in the relative bordism group such that $J[V]$ is
`holomorphic' in $w$. The homomorphism $J$ is constructed using an 
inverse system of truncated projective spaces; around the same time
Mahowald [1,11] called attention to the remarkable properties 
of a similar construction in {\bf framed} bordism, and went on to
consider the properties of this inverse limit \dots \medskip
 
v) The class of the interval $[-1,+1]$ (with sign reversal as
involution) defines an element of the relative bordism group
which maps under $J$ to $w^{-1}$, cf.\ [5]. \bigskip

\noindent
In modern terminology the diagram displays the Tate
$\Z/2\Z$-cohomology [[9]; Swan probably also deserves mention here, cf.\ [14]]
of the forgetful map from geometric $\Z/2\Z$-equivariant unoriented 
cobordism to ordinary unoriented cobordism; this is very closely related 
to the construction used by Kriz [10 Corollary 1.4]
to compute the homotopy-theoretic $\Z/p\Z$-equivariant bordism
groups $MU_{\Z/p\Z*}^{hot}$. In fact the direct limit of the
system defined by suspending the diagram with respect to a
complete family of $\Z/2\Z$-representations (cf.\ [15]) has top
row $$0 \rightarrow \N_{\Z/2\Z*}^{hot} \rightarrow
\N_{*}(BO)[u,u^{-1}] \rightarrow \widetilde \N_{*}(B\Z/2\Z) \rightarrow 
0 \;,$$  which is identical [aside from obvious modifications]
with Kriz's. It follows, in particular, that the
stabilization map $$\N_{\Z/2\Z*}^{geo} \rightarrow
\N_{\Z/2\Z*}^{hot}$$ is injective. 

\section{Complex circle actions}

\noindent
Boardman has an explicit formula for his $J$-homomorphism: if $[V]$ is the 
cobordism class of a closed manifold with involution, then $$J[V] = 
\pi^{-1} \sum_{k \geq 0} w^{k} p_{*} \partial (w^{-k-1} V) \; ,$$ where 
$p_{*} : \N_{*}(B\Z/2\Z) \rightarrow \N_{*}$ sends a manifold with free 
involution to its quotient, and $$\pi := \sum_{k \geq 0} [\R P_{k}]w^{k} \; ,$$
by [4 Theorem 14]. This formula has a natural interpretation in 
terms of formal group laws, but the idea is easier to explain in the context 
of the $\T$-equivariant Tate cohomology of $MU$, which fits in an exact
sequence $$0 \rightarrow MU^{-*}(B\T) \rightarrow t_{\T}(MU_{*})
= MU_{*}((c)) \rightarrow MU_{*-1}(B\T) \rightarrow 0 \; ;$$
here the Chern class $c$ plays the role of the Stiefel-Whitney
class $w$ in the unoriented case. Since the complex cobordism
ring is torsion-free, we can introduce denominators with
impunity: the analogue of $w^{-1}$ is the class $c^{-1}$ of the
two-disk, viewed as a complex-oriented manifold with circle
action free on the boundary. Boardman observes that $$p_{*}
\partial w^{-k-1}  =  [\R P_{k}] \; ;$$ similarly, we have $$p_{*}
\partial c^{-k-1}  =  [\C P_{k}]$$ in the complex-oriented
situation. \bigskip

\noindent {\bf Proposition:} The homomorphism $p_{*} \partial$ is
the formal residue at the origin with respect to an additive
coordinate for the formal group law of $MU$. \medskip

\noindent{\bf Proof:} The idea is to write $c$ as a formal power
series $$c = z \; + \; {\rm terms \; of \; higher \; order} \; ,$$ 
with coefficients chosen so that $$\res_{z=0} c^{-k-1} 
= [\C P_{k}] \; .$$ Clearly $$\sum_{k \geq 0} \res_{0} c^{-k-1} 
u^{k}  =  \res_{0} \frac{c^{-1}}{1-c^{-1}u}  =  \sum_{k \geq
0} [\C P_{k}]u^{k}  =  \log_{MU}'(u) \; .$$ On the other hand, if
$c = f(z)$ then $$\res_{0}(f(z) - u)^{-1} = \frac{1}{2\pi i} \int
\frac{dz}{f(z) - u} \; $$ which equals $$\frac{1}{2\pi i} \int
\frac {(f^{-1})'(c) dc}{c - u} = (f^{-1})'(u) \; ;$$  thus $c =
\exp_{MU}(z)$. \medskip

\noindent {\sc Remarks:} \medskip

i) In the context of remark iv) above, Boardman's formula for $J$
on `holomorphic' elements $V$ is thus essentially the same as
Cauchy's. His use of the symbol $\pi$ suggests that this analogy 
was not far from his mind. \medskip

ii) The formal residue was introduced in Quillen's Bulletin
announcement [12], but seems to have since disappeared from
algebraic topology. Perhaps that paper deserves another look as
well. \medskip

iii) These constructions suggest that while the Chern class $c$
is a natural uniformizing parameter for algebraic questions about
complex cobordism, its inverse may be more natural for geometric
questions. We thus have two reasonable coordinates on cobordism, 
centered at the south and north poles of the Riemann sphere, 
overlapping in the temperate region defined by Tate cohomology. \bigskip

\noindent {\sc Acknowledgements:} I would like to thank John Greenlees, Igor 
Kriz, and Dev Sinha for raising my consciousness about equivariant cobordism,
and Bob Stong for helpful conversations about the history of this subject. \bigskip

\noindent {\sc Postscript:} This has appeared in the Boardman Festschrift\medskip

\noindent
{\bf Homotopy invariant algebraic structures} 15 - 18, Contemp. Math. 239, 
AMS (1999). \medskip

\noindent
I am posting it here in hopes of advertising this geometric description of
Tate cohomology. Similar ideas play a role in\medskip

\noindent
Heisenberg groups and algebraic topology, in the Segal Festscrift {\bf 
Topology, geometry and quantum field theory} 235 - 246, LMS Lecture Notes 
308, Cambridge (2004); as well as in \medskip

\noindent
Completions of $\Z/(p)$-Tate cohomology of periodic spectra. Geom. Topol. 2 (1998) 145 - 174,
{\tt arXiv:math/9808141} 

\newpage

\bibliographystyle{amsplain}

\end{document}